\documentclass[letterpaper, 10 pt, journal, twoside]{ieeetran}

\IEEEoverridecommandlockouts                              

\usepackage{amsmath} 	
\usepackage{amssymb}  	
\usepackage{amsfonts}
\usepackage{graphicx}
\usepackage{subfigure}
\usepackage{epstopdf}
\usepackage{dblfloatfix}
\usepackage{tikz}
\usetikzlibrary{shapes,matrix,decorations.shapes}
\usepackage{array,xcolor}
\usepackage{algorithm}
\usepackage{algpseudocode}
\usepackage{multirow}
\usepackage{threeparttable}
\usepackage{booktabs}
\usepackage{balance}
\usepackage{cite}

\usepackage[colorlinks=true,      
			linkcolor=black,      
			citecolor=black,      
			filecolor=black,      
			urlcolor=blue ]{hyperref}

\newtheorem{definition}{Definition}

\newtheorem{proposition}{Proposition}
\newtheorem{remark}{Remark}
\newtheorem{lemma}{Lemma}

\begin{document}

\title{Exploiting Sparsity in the Coefficient Matching Conditions in Sum-of-Squares Programming \\ using ADMM}


\author{Yang Zheng,~Giovanni~Fantuzzi and~Antonis Papachristodoulou,~\IEEEmembership{Senior Member,~IEEE} 
\thanks{Y. Zheng is supported by the Clarendon Scholarship and the Jason Hu Scholarship. A. Papachristodoulou is supported by project EP/M002454/1.}
\thanks{Y. Zheng and A. Papachristodoulou are with the Department of Engineering Science, University of Oxford, Parks Road, Oxford, OX1 3PJ, U.K. (e-mail: \{yang.zheng, antonis\}@eng.ox.ac.uk)}%
\thanks{G.~Fantuzzi is with the Department of Aeronautics, Imperial College London, South Kensington Campus, London, SW7 2AZ, United Kingdom (e-mail: gf910@ic.ac.uk)}%
}

\maketitle
\pagenumbering{gobble}
\begin{abstract}
    This paper introduces an efficient first-order method based on the alternating direction method of multipliers (ADMM) to solve semidefinite programs (SDPs) arising from sum-of-squares (SOS) programming. We exploit the sparsity of the \emph{coefficient matching conditions} when SOS programs are formulated in the usual monomial basis to reduce the computational cost of the ADMM algorithm. Each iteration of our algorithm requires one projection onto the positive semidefinite cone and the solution of multiple quadratic programs with closed-form solutions free of any matrix inversion. Our techniques are implemented in the open-source MATLAB solver SOSADMM. Numerical experiments on SOS problems arising from unconstrained polynomial minimization and from Lyapunov stability analysis for polynomial systems show speed-ups compared to the interior-point solver SeDuMi, and the first-order solver CDCS.
\end{abstract}


\section{Introduction}

\IEEEPARstart{C}{hecking} whether a given polynomial is nonnegative has applications in many areas (see~\cite{lasserre2009moments,chesi2009guest} and the references therein). For example, the unconstrained polynomial optimization problem $\min_{x\in\mathbb{R}^n} p(x)$ is equivalent to 
\begin{equation} \label{eq:POP}
    \begin{aligned}
        \max & \quad \gamma \\
        \text{subject to} & \quad p(x)- \gamma \geq 0.
    \end{aligned}
\end{equation}
Moreover, the stability of an equilibrium $x^*$ of a polynomial dynamical system $\dot{x}(t) = f(x(t))$, $x(t)\in\mathbb{R}^n$ in a neighbourhood $\mathcal{D}$ of {$x^*$} (we assume $x^*=0$  without loss of generality)---a fundamental problem in control theory---can be established by constructing a polynomial $V(x)$ (called {\it Lyapunov function}) that satisfies the polynomial inequalities~\cite{papachristodoulou2005tutorial}
\begin{equation} \label{eq:FindLyapunov}
    \begin{cases}
        \phantom{-}V(x) > 0, &\forall x \in \mathcal{D} \backslash \{0\}, \\
        -\dot{V}(x) = - \langle \nabla V(x), f(x) \rangle \geq 0, &\forall x \in \mathcal{D}.
    \end{cases}
\end{equation}
Throughout this work, $\langle \cdot,\cdot \rangle$ denotes the inner product in the appropriate Hilbert space.

A powerful way to test polynomial inequalities is to employ a sum-of-squares (SOS) relaxation (we refer the reader to~\cite{parrilo2003semidefinite,lasserre2001global} for details on SOS relaxations in polynomial optimization, and to~\cite{papachristodoulou2005tutorial} for a tutorial on SOS techniques for systems analysis). In fact, while testing the non-negativity of a polynomial is NP-hard in general, the existence of an SOS decomposition can be checked in polynomial time by solving a semidefinite program (SDP)~\cite{parrilo2003semidefinite}. Unfortunately, however, the size of the SDP for the SOS relaxation of a degree-$d$ polynomial in $n$ variables is $\begin{pmatrix}\begin{smallmatrix}\! n+d \\ d \\ \end{smallmatrix} \end{pmatrix}$. Consequently, SOS relaxations are limited to small problem instances; with the current technology, for example, Lyapunov-based analysis is impractical for general systems with ten or more states.

In order to mitigate scalability issues, one can act at the \emph{modeling} level, \emph{i.e.} one can try to replace the SDP obtained from an SOS relaxation with an optimization problem that is cheaper to solve still using second-order interior-point methods (IPMs), implemented in efficient solvers such as SeDuMi~\cite{sturm1999using}. One approach is to exploit structural properties of the polynomial whose positivity is being tested~\cite{reznick1978extremal,lofberg2009pre, gatermann2004symmetry, waki2006sums,permenter2014basis}. For example, computing the Newton polytope~\cite{reznick1978extremal} or checking for diagonal inconsistency~\cite{lofberg2009pre} can restrict the monomial basis required in the SOS decomposition by eliminating redundant monomials. Further improvements are possible by group-theoretic symmetry reduction techniques~\cite{gatermann2004symmetry} and graph-theoretic correlative sparsity~\cite{waki2006sums}. Facial reduction has also been applied to select a reduced monomial basis for SOS programs in~\cite{permenter2014basis}. A second approach is to approximate the positive semidefinite (PSD) cone using diagonally dominant or scaled diagonally dominant matrices~\cite{majumdar2014control, ahmadi2014dsos}. These relaxations can be solved with linear programs (LPs) or second-order-cone programs (SOCPs), rather than SDPs, and the conservativeness introduced by approximating the PSD cone can be reduced with a recently proposed basis pursuit algorithm~\cite{ahmadi2015sum}.

Further improvements are available on the \emph{computational} level if IPMs are replaced by more scalable first-order methods (FOMs) at the cost of reduced accuracy. The design of efficient first-order algorithms for large-scale SDPs has received particular attention in recent years. For instance, Wen \emph{et al.} proposed an alternating direction augmented Lagrangian method for large-scale dual SDPs~\cite{wen2010alternating}. O'Donoghue \emph{et al.} developed an operator-splitting method to solve the homogeneous self-dual embedding of conic programs~\cite{ODonoghue2016}, which has recently been extended by the authors to exploit aggregate sparsity via chordal decomposition~\cite{ZFPGW2016,zheng2016fast}. In the context of SOS programming, Bertsimas \emph{et al.} proposed an accelerated FOM for unconstrained polynomial optimization~\cite{bertsimas2013accelerated}, while Henrion \& Malick introduced a projection-based method for SOS relaxations~\cite{henrion2012projection}. However, both approaches are only applicable to a small subset of SOS programs as they rely on the constraint of the corresponding SDP being orthogonal, which is not the case for SOS problems with free variables.

In this paper, we propose a first-order algorithm based on the alternating direction method of multipliers (ADMM) to solve the SDPs arising from SOS optimization. In contrast to~\cite{bertsimas2013accelerated} and~\cite{henrion2012projection}, we exploit the sparsity in the coefficient matching conditions, making our approach suitable for a larger class of SOS programs. While the aggregate sparsity pattern of these SDPs is dense (so that the methods of~\cite{ZFPGW2016,zheng2016fast} are not very advantageous), each equality constraint in the SDP only involves a small subset of decision variables when an SOS program is formulated in the usual monomial basis. This sparsity can be exploited to formulate an efficient ADMM algorithm, the iterations of which consist of conic projections and optimization problems with closed-form solutions that---crucially---are free of any matrix inversion. We implement our techniques in SOSADMM, an open-source MATLAB solver. The efficiency of our methods compared to the IPM solver SeDuMi~\cite{sturm1999using} and the first-order solver CDCS~\cite{CDCS} is demonstrated on SOS problems arising from unconstrained polynomial optimization and from Lyapunov stability analysis of polynomial systems.

The rest of this paper is organized as follows. Section~\ref{se:preliminaries} reviews SOS polynomials and the ADMM algorithm. Sparsity for SDPs arising in SOS programs is discussed in Section~\ref{se:sparsity}, and we show how to exploit it to build an efficient ADMM algorithm in Section~\ref{se:ADMM}. Numerical experiments are reported in Section~\ref{se:simulation}. Section~\ref{se:conclusion} concludes the paper.

\section{Preliminaries}\label{se:preliminaries}


\subsection{SOS polynomials and SDPs}

The sets of real and natural numbers {(including zero)} are denoted by $\mathbb{R}$ and $\mathbb{N}$, respectively. Let $x \in \mathbb{R}^n$, $\alpha \in \mathbb{N}^n$, and let $ x^{\alpha}= x_1^{\alpha_1}x_2^{\alpha_2}\cdots x_n^{\alpha_n} $ denote a monomial in $ x $ of degree $ \lvert \alpha \rvert = \sum_{i = 1}^n \alpha_i $. Given an integer $ d \in \mathbb{N} $, we denote $ \mathbb{N}_d^n = \{ \alpha \in \mathbb{N}^n : \vert \alpha\vert  \leq d \} $, and the vector of all monomials of degree no greater than $ d $ by
\begin{equation}\label{eq:CandidateMonomials}
    \begin{aligned}
        v_d(x) &= \{x^{\alpha} \mid \alpha \in \mathbb{N}_d^n\} \\
                & =[ 1,x_1,x_2,\ldots,x_n,x_1^2,x_1x_2,\ldots,x_n^d ]^T.
  \end{aligned}
\end{equation}
The length of $ v_d(x) $ is $ |\mathbb{N}_d^n|=\begin{pmatrix}\begin{smallmatrix} \! n+d \\ d \\ \end{smallmatrix} \end{pmatrix} $.
A real polynomial $p(x)$ is a finite, real linear combination of monomials of $x$
\begin{equation*}
  p(x) = \sum_{\alpha \in \mathbb{N}^n} p_{\alpha}x^{\alpha},
  \qquad p_{\alpha} \in \mathbb{R}.
\end{equation*}
The degree of $p(x)$ is the maximum of the degrees of all monomials with nonzero coefficients. We  denote the set of real polynomials in $x$ by $\mathbb{R}[x]$. 
\begin{definition}
    A polynomial $p(x)\in\mathbb{R}[x]$ of degree $2d$ is a sum-of-squares (SOS) if there exist polynomials $ f_i(x) \in \mathbb{R}[x]$, $i = 1,\ldots, m $ of degree no greater than $d$ such that
    $$
        p(x) = \sum_{i = 1}^m \left[f_i(x)\right]^2.
    $$
\end{definition}
\vspace{0.5em}
Clearly, the existence of an SOS representation guarantees that $ p(x)\geq 0 $. The following theorem gives an equivalent characterization of SOS polynomials.
\begin{proposition}[$\!\!$\cite{parrilo2003semidefinite}]
\label{Pro:SOS}
    A polynomial $ p(x) \in \mathbb{R}[x] $ of degree $2d$ is an SOS polynomial if and only if there exists a $\begin{pmatrix}\begin{smallmatrix} \! n+d \\ d \\ \end{smallmatrix} \end{pmatrix} \times  \begin{pmatrix}\begin{smallmatrix} \!  n+d \\ d \\ \end{smallmatrix} \end{pmatrix} $ symmetric PSD matrix $ X \succeq 0$ such that
    \begin{equation} \label{eq:SOSform}
        p(x) = v_d(x)^TXv_d(x).
    \end{equation}
\end{proposition}
    \vspace{6pt}

The equality in~\eqref{eq:SOSform} gives a set of affine equalities on the elements of $X$ to match the coefficients of $p(x)$. Together with $X\succeq 0$, this makes the problem of finding an SOS representation for $p(x)$ an SDP. The formulation of such SDPs can be assisted by software packages, such as SOSTOOLS \cite{papachristodoulou2013sostools} and GloptiPoly~\cite{henrion2003gloptipoly}.

\begin{remark}
    The size of the PSD matrix $X$ in~\eqref{eq:SOSform} is $ \begin{pmatrix}\begin{smallmatrix}\!  n+d \\ d \\ \end{smallmatrix} \end{pmatrix} \times  \begin{pmatrix}\begin{smallmatrix} \! n+d \\ d \\ \end{smallmatrix} \end{pmatrix} $ because we have used the full set of monomials of degree no greater than $ d $ in our representation. This number might be reduced by inspecting the structural properties of $p(x)$ to identify and eliminate redundant monomials in $v_d(x)$; well-known techniques include Newton polytope~\cite{reznick1978extremal}, diagonal inconsistency~\cite{lofberg2009pre}, symmetry property~\cite{gatermann2004symmetry}, and facial reduction~\cite{permenter2014basis}.
\end{remark}

\begin{table*}
    \centering
    \setlength{\abovecaptionskip}{0pt}
    \setlength{\belowcaptionskip}{1em}
    \renewcommand\arraystretch{0.8}
    \caption{Density of nonzero elements in the equality constraints of SDP~\eqref{eq:SDPSOS}}
    \label{T:density}
    \begin{tabular}{c| c c c c c c c}
        \hline \toprule[1pt] 
         $n$   & 4 & 6  & 8 & 10 & 12 & 14 & 16\\
        \hline\\[-0.75em]
        $2d = 4$ & $1.42\times 10^{-2}$ & $4.76\times 10^{-3}$ & $2.02\times 10^{-3}$ & $9.99\times 10^{-4}$  & $5.49\times 10^{-4}$& $3.27\times 10^{-4}$ & $2.06\times 10^{-4}$\\
        $2d = 6$ & $4.76\times 10^{-3}$ & $1.08\times 10^{-3}$ & $3.33\times 10^{-4}$ & $1.25\times 10^{-4}$  & $5.39\times 10^{-5}$& $2.58\times 10^{-5}$ & $1.34\times 10^{-5}$\\
        $2d = 8$ & $2.02\times 10^{-3}$ & $3.33\times 10^{-4}$ & $7.77\times 10^{-5}$ & $2.29\times 10^{-5}$  & $7.94\times 10^{-6}$ & $3.13\times 10^{-6}$ & $1.36\times 10^{-6}$\\
        \bottomrule[1pt]
        \end{tabular}
\end{table*}

\subsection{ADMM algorithm}

The ADMM algorithm solves the optimization problem
\begin{equation}
\label{e:ADMM}
    \begin{aligned}
        \min_{y,z} \quad & f(y)+g(z) \\
        \text{subject to} \quad & Ay + Bz = c,
    \end{aligned}
\end{equation}
where $y\in\mathbb{R}^n$ and $z\in\mathbb{R}^m$ are the decision variables, $f:\mathbb{R}^n\to\mathbb{R}$ and $g:\mathbb{R}^m\to\mathbb{R}$ are convex functions, and $A\in\mathbb{R}^{l\times n}$, $B\in\mathbb{R}^{l\times m}$ and $c\in\mathbb{R}^{l}$ are the constraint data. Given a penalty parameter $\rho>0$ and a multiplier $\lambda\in\mathbb{R}^{l}$ (known as the \emph{dual variable}), the ADMM algorithm
solves~\eqref{e:ADMM} by finding a saddle point{~\cite[Chapter 5.4]{boyd2004convex}} of the augmented Lagrangian
\begin{multline} \label{E:AugLag}
    L_{\rho}(y,z,\lambda)  = f(y) + g(z) + \lambda^T\, (Ay + Bz -c) \\
                    + \frac{\rho}{2} \left\|Ay + Bz - c\right\|^2
\end{multline}
with the following steps:
\begin{subequations}\label{E:ADMM}
    \begin{align}
        y^{k+1} & = \text{arg} \min_{y} L_{\rho}(y,z^{k},\lambda^{k}),
        \label{E:ADMM_S1}\\
        z^{k+1} & = \text{arg} \min_{z} L_{\rho}(y^{k+1},z,\lambda^{k}),
        \label{E:ADMM_S2}\\
        \lambda^{k+1} &= \lambda^{k} + \rho ( A y^{k+1} + B z^{k+1} - c).\label{E:ADMM_S3}
    \end{align}
\end{subequations}
In these equations, the superscript $k$ denotes the value of a variable at the $k$-th iteration of the algorithm, and $\|\cdot\|$ denotes the standard Euclidean norm, \emph{i.e.}, $\|x\| = \sqrt{x^Tx}$ for $x \in \mathbb{R}^n$. Then, from a computational perspective, steps~\eqref{E:ADMM_S1} and~\eqref{E:ADMM_S2} are equivalent to the minimizations of
    \begin{equation}\label{E:scaleAugLagPaper}
         \tilde{L}_{\rho}(y,z,\lambda) = f(y) + g(z)
         + \frac{\rho}{2} \left\|Ay + Bz - c + \frac{1}{\rho} \lambda\right\|^2
    \end{equation}
    over $y$ and $z$, respectively, with $\lambda$ fixed.
Under very mild conditions, ADMM converges to a solution with a rate $\mathcal{O}(\frac{1}{k})$, which is independent of $\rho$, although its value can affect convergence in practice; see~\cite[Section 3.2]{boyd2011distributed} for details.

\section{Row Sparsity in SDPs from SOS programs} \label{se:sparsity}

\subsection{SDP formulations of SOS relaxations}
Let $A_\alpha $ be the indicator matrix for the monomials $x^{\alpha}$ in the rank-one matrix $v_d(x)v_d(x)^T$; in other words, the entry of $A_\alpha $ with row index $\beta$ and column index $\gamma$ (where the natural ordering for multi-indices $\beta,\gamma\in\mathbb{N}_d^n$ is used) satisfies
\begin{equation} \label{eq:coeff}
    (A_\alpha)_{\beta,\gamma} =
    \begin{cases}
    		1 &\text{if } \beta + \gamma = \alpha\\
    		0 &\text{otherwise}.
    	\end{cases}
\end{equation}
The SOS constraint \eqref{eq:SOSform} can then be reformulated as
\begin{equation} 
    p(x) = \langle v_d(x)v_d(x)^T, X \rangle
    		 = \sum_{\alpha\in\mathbb{N}_{2d}^n} \langle A_\alpha, X\rangle x^\alpha.
\end{equation}
Matching the coefficients of the left- and right-hand sides gives the equality constraints
\begin{equation} \label{eq:CoeffCondition1}
        \langle A_\alpha,X \rangle = p_{\alpha} \quad \forall \: \alpha \in \mathbb{N}^n_{2d}.
\end{equation}
We refer to these equalities as \textit{coefficient matching conditions}. The existence of an SOS decomposition for $p(x)$ (or lack thereof) can then be checked with the feasibility SDP
\begin{equation}
\label{eq:SDPSOS}
    \begin{aligned}
        \text{find}\quad &X \\
        \text{subject to} \quad  & \langle A_{\alpha}, X \rangle = p_{\alpha},
        							\quad \alpha\in\mathbb{N}_{2d}^n,\\
            & X \succeq 0.
    \end{aligned}
\end{equation}
When the full monomial basis is used, as in this case, the dimension of $X$ and the number of constraints in~\eqref{eq:SDPSOS} are, respectively,
\begin{align}
\label{eq:SDPdimension}
N = |\mathbb{N}^n_{d}| &= \begin{pmatrix} n+d \\d\end{pmatrix},
&
m = |\mathbb{N}^n_{2d}| &= \begin{pmatrix} n+2d \\2d\end{pmatrix}.
\end{align}

\subsection{Properties of the coefficient matching conditions}
\label{se:sparCoefficient}
In this section, for simplicity, we re-index the constraint matching conditions~\eqref{eq:CoeffCondition1} using integer indices $i=1,\ldots,m$ instead of the multi-indices $\alpha$.

The conditions~\eqref{eq:CoeffCondition1} inherit two important properties from the data matrices $A_{i}$, $i=1,\ldots,m$. The first one follows from the fact that the matrices $A_{i}$ are orthogonal. If $n_i$ denotes the number of nonzero entries in $A_i$ we have
\begin{equation} \label{eq:coeffOrthogonal}
    \langle A_i, A_{j}\rangle =
    \begin{cases}
    		n_{i} &\text{if } i = j,\\
    	    0 &\text{otherwise}.
    	\end{cases}
\end{equation}
After letting $\mathrm{vec}: \mathbb{S}^N \to \mathbb{R}^{N^2}$ be the usual operator mapping a matrix to the stack of its columns, and defining
\begin{align} \label{E:VecMatrixDef}
    A = \begin{bmatrix} \text{vec}(A_1) & \cdots & \text{vec}(A_m) \end{bmatrix}^T,
\end{align}
the equality constraints in~\eqref{eq:SDPSOS} can be rewritten as the matrix-vector product $A\cdot \text{vec}(X) = b$,
where $b\in \mathbb{R}^m$ is a vector collecting the coefficients $p_{i}$, $i=1,\ldots,m$. Property~\eqref{eq:coeffOrthogonal} directly implies the following lemma, which forms the basis of the FOMs of~\cite{bertsimas2013accelerated,henrion2012projection}.
\begin{lemma}[Orthogonality of constraints] \label{lemma:orthogonality}
    $AA^T$ is an $m\times m$ diagonal matrix with $(AA^T)_{ii} = n_i$.
\end{lemma}

The second property of the coefficient matching conditions is that they are sparse, in the sense that each equality constraint in~\eqref{eq:SDPSOS} only involves a small subset of entries of $X$, because only a small subset of entries of the product $v_d(x)v_d(x)^T$ are equal to a given monomial $x^\alpha$. Thus, the vectorized matrix $A$ is \emph{row sparse}, meaning that each row is a sparse vector. In particular, {the following result holds}.
\begin{lemma}[Sparsity of constraints] \label{lemma:sparsity}
    Let $A$ be the vectorized matrix for~\eqref{eq:SDPSOS}, and let $N$ and $m$ be as in~\eqref{eq:SDPdimension}. The number of nonzero elements in $A$ is $N^2$, and the density of nonzero elements in $A$ is equal to
$m^{-1} = \mathcal{O}({n^{-2d}})$.
\end{lemma}

\begin{IEEEproof}
Since the matrix  $v_d(x)v_d(x)^T$ contains all monomials $x^\alpha$, $\alpha\in\mathbb{N}_{2d}^n$, all entries of the PSD matrix $X$ enter at least one of the equality constraints in~\eqref{eq:SDPSOS}. Moreover,~\eqref{eq:coeffOrthogonal} implies that each entry of $X$ enters at most one constraint. Therefore, $A$ must contain $N^2$ nonzero elements. Its density is then given by
\begin{equation*}
\frac{N^2}{N^2\times m} = \frac{1}{m} = \left[\begin{pmatrix} n+2d \\2d \end{pmatrix}\right]^{-1} = \mathcal{O}({n^{-2d}}).
\end{equation*}
\end{IEEEproof}

\begin{figure}[t]
	\setlength{\abovecaptionskip}{0pt}
    \setlength{\belowcaptionskip}{0em}
    \centering
    { \label{fig:a}
	\includegraphics[scale=1.125]{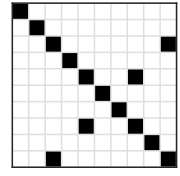}
    }
    \caption{Sparsity pattern of $AA^T$ for the example~\eqref{eq:ExNonOrth}}
    \label{fig:NonOrtho}
\end{figure}

\begin{remark} \label{r:remark2}
While the constraint matrix $A$ is sparse (see typical values in Table~\ref{T:density}), the aggregate sparsity pattern of the SDP~\eqref{eq:SDPSOS} is dense because all entries of the matrix variable $X$ appear in the equality constraints. This implies that $X$ is generally a dense variable, so the first-order algorithms of~\cite{ZFPGW2016,zheng2016fast} are not particularly suitable.
\end{remark}

\begin{remark}\label{r:remark3}
    The property of orthogonality in Lemma~\ref{lemma:orthogonality} holds for standard SOS feasibility problems. However, this property fails for the following example:
    \begin{equation} \label{eq:ExNonOrth}
        \begin{aligned}
            \text{find} & \quad a,\,b \\
            \text{subject to} & \quad ax^4 + bx^2 + x + 1 \; \text{is SOS},\\
                             & \quad bx^4 + ax^2 + x + 1 \; \text{is SOS}.
        \end{aligned}
    \end{equation}
     Fig.~\ref{fig:NonOrtho} shows the sparsity pattern of $AA^T$ for~\eqref{eq:ExNonOrth} obtained using SOSTOOLS, demonstrating that the constraints are not orthogonal. The reason is that~\eqref{eq:ExNonOrth} involves the free parameters $a,\,b$ as well as the PSD matrices for the SOS representation, and this destroys the orthogonality. This issue is common in control applications; see, \emph{e.g.}, the condition~\eqref{eq:FindLyapunov} when finding Lyapunov functions. Consequently, the first-order algorithms in~\cite{bertsimas2013accelerated,henrion2012projection} cannot be applied to many problems with SOS constraints because they rely on constraint orthogonality.
\end{remark}

\section{Exploiting Row Sparsity in SDPs} \label{se:ADMM}

As we have seen, the algorithms in~\cite{ZFPGW2016,zheng2016fast,bertsimas2013accelerated,henrion2012projection} are not useful for generic SOS problems because their aggregate sparsity pattern is dense, and the orthogonality property only holds for simple SOS feasibility problems. However, the data matrix $A$ is always row-sparse due to the coefficient matching conditions in SOS programs. This property can be exploited to construct an efficient ADMM algorithm. 
In the following, we consider a generic SDP in the vectorized form
\begin{equation} \label{E:PrimalVectorForm}
\begin{aligned}
    \min_{x} \quad & c^Tx\\
    \text{subject to} \quad & Ax=b, &
					   & x \in \mathcal{K},
\end{aligned}
\end{equation}
where $x$ is the optimization variable, $A\in \mathbb{R}^{m \times \hat{n}}$, $b \in \mathbb{R}^{m}$ and $c\in \mathbb{R}^{\hat{n}}$ are the problem data, and $\mathcal{K}$ is a product of cones, at least one of which is the PSD cone.  Throughout this section, $\delta_{\mathcal{S}}(x)$ denotes the indicator function of a set $\mathcal{S}$,
$$
    \delta_{\mathcal{S}}(x) =
    \begin{cases}
    0, &\text{if } x \in \mathcal{S}, \\
    + \infty, &{\text{if } x \notin \mathcal{S}}.
    \end{cases}
$$

\subsection{Reformulation considering individual row sparsity}

Let us represent $A = [a_1, a_2, \ldots, a_m]^T$, so each vector $a_i$ is a row of $A$, and let $H_i \in \mathbb{R}^{|a_i|\times \hat{n}}, i=1,\ldots,m$ be ``entry-selector''
matrices of 1's and 0's selecting the nonzero elements of $a_i$, where $|a_i|$ denotes the number of nonzero elements of $a_i$. Note that the rows of $H_i$ are orthonormal, since each selects a different entry of $a_i$. Then,
\begin{equation} \label{eq:LocalVariable}
    Ax = b \Leftrightarrow
    \begin{cases}
    (H_ia_i)^Tz_i = b_i, &i = 1, \ldots, m, \\
    z_i = H_ix,&i = 1, \ldots, m.
    \end{cases}
\end{equation}
In \eqref{eq:LocalVariable}, $z_i$ is a copy of the elements of $x$ which enter the $i$-th affine constraint. 
It is also convenient to introduce an additional slack variable $u=x$, so the affine constraints in~\eqref{eq:LocalVariable} and conic constraint in~\eqref{E:PrimalVectorForm} are decoupled when applying the ADMM algorithm. We can then reformulate~\eqref{E:PrimalVectorForm} as
\begin{equation} \label{eq:RowDecomposedForm}
\begin{aligned}
    \min_{z_i,u,x} \quad & c^Tx\\
    \text{subject to} \quad & (H_ia_i)^Tz_i = b_i & i = 1, \ldots, m,\\
                            & z_i = H_ix, & i = 1, \ldots, m,\\
                            & u = x,
					         u \in \mathcal{K}.
    \end{aligned}
\end{equation}

\subsection{ADMM steps}

To apply ADMM, we move the affine constraints $(H_ia_i)^Tz_i = b_i$ and the conic constraint $u \in \mathcal{K}$ in \eqref{eq:RowDecomposedForm} to the objective using the indicator functions $\delta_0(\cdot)$ and $\delta_{\mathcal{K}}(\cdot)$, respectively:
\begin{equation} \label{eq:ADMMPrimal}
    \begin{aligned}
    \min_{z_i,u,x} \quad &c^Tx + \delta_{\mathcal{K}}(u) + \sum_{i=1}^{m} \delta_{0}\left((H_ia_i)^Tz_i - b_i\right) \\
    \text{subject to} \quad & z_i = H_ix, \quad i = 1, \ldots, m,\\
    & u = x.
    \end{aligned}
\end{equation}

The augmented Lagrangian of \eqref{eq:ADMMPrimal} is
\begin{multline}
    \label{eq:AugLag1}
        L =
        c^Tx
        + \delta_{\mathcal{K}}(u)
        + \sum_{i=1}^m\delta_0\!\left[(H_i a_i)^Tz_i - b_i\right]
        \\
        + \sum_{i=1}^m
        \mu_i^T\left( z_i - H_ix \right)
        +\frac{\rho}{2}\left\|z_i - H_ix\right\|^2
        \\
        + \xi^T\left( u-x \right)
        + \frac{\rho}{2}\left\|u -x\right\|^2,
\end{multline}
and we group the variables as
\begin{equation*}
    \begin{aligned}
        \mathcal{Y} = \{x\}, \;
        \mathcal{Z} = \{u,z_1,\ldots,z_m\}, \;
        \mathcal{D} = \{\mu_1,\ldots,\mu_m,\xi\}.
    \end{aligned}
\end{equation*}
Then, the ADMM steps~\eqref{E:ADMM_S1}--\eqref{E:ADMM_S3} become the following:

\subsubsection{Minimization over $\mathcal{Y}$} The minimization of \eqref{eq:AugLag1} over the variables in $\mathcal{Y}$ is {equivalent to} an unconstrained quadratic program,
\begin{equation}
    \label{eq:minBlockX1}
        \min_{x}  c^Tx
        +\frac{\rho}{2}\sum_{i=1}^m\!
        \left\|z_i^k\!-\!H_i x\!+\!\frac{\mu_i^k}{\rho}\right\|^2
        \!+\!\frac{\rho}{2}\left\|u^k\!-\!x\!+\!\frac{\xi^k}{\rho}\right\|^2\!\!.
\end{equation}
The updated variable $x^{k+1}$ is then simply given by
\begin{equation} \label{eq:ADMMBlkX1}
    \!x^{k+1}\!=\!D^{-1}\!\left[
    \sum_{i=1}^m\!H_i^T\!\left(z_i^k\!+\!\frac{\mu_i^k}{\rho}\right)
    \!+\!\left(u^k\!+\!\frac{\xi^k}{\rho}\right)
    -\!\frac{1}{\rho}c \right]\!,
\end{equation}
where the matrix
$
    D = I+\sum_{i=1}^m H_i^TH_i
$
is diagonal because the rows of each matrix $H_i$ are orthonormal. This means that~\eqref{eq:ADMMBlkX1} is cheap to calculate.
\subsubsection{Minimization over $\mathcal{Z}$} Minimizing \eqref{eq:AugLag1} over the variables in $\mathcal{Z}$ amounts to a conic projection,
\begin{equation}
    \label{eq:minBlockYz1}
    \begin{aligned}
        \min_{u}  & \quad \left\|u - x^{k+1} + {\rho^{-1}\xi^{k}}\right\|^2 \\
        \text{subject to} & \quad  u \in \mathcal{K},
    \end{aligned}
\end{equation}
plus $m$ independent quadratic programs
\begin{equation}
    \label{eq:minBlockYzi1}
    \begin{aligned}
        \min_{z_i} \quad  &\left\|z_i - H_ix^{k+1} + {\rho^{-1}\mu_i^k}\right\|^2 \\
        \text{subject to} \quad & (H_i a_i)^Tz_i = b_i.
    \end{aligned}
\end{equation}
The projection \eqref{eq:minBlockYz1} is easy to compute when $\mathcal{K}$ is a product of $\mathbb{R}^n$, the non-negative orthant, second-order cones, and PSD cones; for example, a projection onto the PSD cone only requires one eigen-decomposition. As for problem~\eqref{eq:minBlockYzi1}, its KKT conditions are
\begin{subequations}
    \begin{align}
        \label{E:KKT1}
        z_i - H_i x^{k+1} + {\rho^{-1} \mu_i^k}+ (H_i a_i)\omega_i &= 0,\\
            \label{E:KKT2}
        (H_i a_i)^Tz_i &= b_i,
    \end{align}
\end{subequations}
where $\omega_i$ is the Lagrangian multiplier for the equality constraint in~\eqref{eq:minBlockYzi1}. Simple algebra shows that
\begin{equation*} 
    \omega_i = 
    \frac{1}{\|H_ia_i\|^{2}}
    \left(- b_i +  (H_i a_i)^TH_i x^{k+1} -
    \frac{1}{\rho}(H_i a_i)^T\mu_i^k
    \right),
\end{equation*}
so the solution $z_i^{k+1}$ to~\eqref{eq:minBlockYzi1} can be calculated easily with~\eqref{E:KKT1}. 
Note that this step is free of any matrix inversion.

\subsubsection{Update multipliers $\mathcal{D}$}
{
According to~\eqref{E:ADMM_S3}, the multipliers in $\mathcal{D}$ are updated with inexpensive and parallelizable gradient ascent steps:
\begin{equation}
    \label{E:Multipliers}
    \begin{aligned}
        \mu_i^{k+1} &= \mu_i^{k} + \rho(z_i^{k+1} - H_ix^{k+1}),
        &i = 1,\ldots,m, \\
        \xi^{k+1} &= \xi^{k} + \rho(u^{k+1} - x^{k+1}).
    \end{aligned}
\end{equation}
}

\subsection{Summary of the computations in the ADMM algorithm}

In the proposed ADMM algorithm, subproblems~\eqref{E:ADMM_S1} and \eqref{E:ADMM_S2} have explicit closed-form solutions. Each iteration requires solving
\begin{enumerate}
  \item one unconstrained quadratic program, given by~\eqref{eq:minBlockX1};
  \item one conic projection, given by~\eqref{eq:minBlockYz1};
  \item $m$ independent quadratic programs, given by~\eqref{eq:minBlockYzi1}.
\end{enumerate}
Note that only the nonzero elements of $a_i$ appear in~\eqref{eq:minBlockYzi1}. Since we have assumed that $a_i$ is sparse, only operations on vectors of small size are required. Besides, our algorithm is free of matrix inversion (with the exception of the $m\times m$ diagonal matrix $D$, requiring $O(m)$ flops), which results from introducing the local variables $z_i$ so each affine constraint can be considered individually. The cost is that our algorithm needs to maintain multiple local variables $z_i$, which may have adverse effects on the convergence speed of ADMM.

In contrast, the FOMs in~\cite{bertsimas2013accelerated, henrion2012projection} fail to deal with general SOS programs since they rely on orthogonality of constraints, and those in~\cite{ZFPGW2016,zheng2016fast} require factorizing the $m \times m$ matrix $AA^T$ (in general, $O(m^3)$ flops). In SDPs arising from generic SOS relaxations, this step is computationally demanding because $AA^T$ is not diagonal (as seen in Remark~\ref{r:remark3}, this is often the case) and the number $m$ is usually large ($m = 18564$ if $n = 12$ and $2d = 6$ in~\eqref{eq:SDPSOS}). Note that all these algorithms converge at rate $\mathcal{O}(\frac{1}{k})$ because they are based on ADMM or its variants.

\section{Numerical Experiments} \label{se:simulation}

We implemented our techniques in SOSADMM, an open-source first-order MATLAB solver for conic programs with row sparsity. Currently, SOSADMM supports cartesian products of the following cones: $\mathbb{R}^n$, non-negative orthant, second-order cone, and the positive definite cone. SOSADMM, the numerical examples presented in this section, and additional examples are available from

\small
\url{https://github.com/oxfordcontrol/SOSADMM}
\normalsize

We tested SOSADMM on random unconstrained polynomial optimization problems and Lyapunov stability analysis of polynomial systems. To assess the suboptimality of the solution returned by SOSADMM, we compared it to the accurate one computed with the interior-point solver SeDuMi~\cite{sturm1999using}. CPU times were compared to the first-order solver CDCS~\cite{CDCS}, which exploits aggregate sparsity in SDPs; in particular, the primal method in CDCS was used~\cite{ZFPGW2016}. In our experiments, the termination tolerance for SOSADMM and CDCS was set to $10^{-4}$, and the maximum number of iterations was $2000$. To improve convergence, SOSADMM employs an adaptive penalty parameter update rule~\cite{boyd2011distributed}, with an initial value $\rho=1$. All tests were run on a PC with a 2.8 GHz Intel\textsuperscript{\textregistered} Core\textsuperscript{\texttrademark} i7 CPU and 8GB of RAM.

\subsection{Unconstrained polynomial optimization}
Consider the global polynomial minimization problem
\begin{equation} \label{eq:ExPOP}
    \min_{x\in\mathbb{R}^n} \; p(x),
\end{equation}
where $p(x)$ is a given polynomial. This problem is equivalent to~\eqref{eq:POP}, and we can obtain an SDP relaxation by replacing the non-negativity constraint with an SOS condition on $p(x) - \gamma $. Motivated by~\cite{henrion2012projection}, we generated $p(x)$ according to
$$
    p(x) = p_0(x) + \sum_{i=1}^n x_i^{2d},
$$
where $p_0(x)$ is a random polynomial of degree strictly less than $2d$. We used GloptiPoly~\cite{henrion2003gloptipoly} to generate the examples.

Table~\ref{T:TimePOP} compares the CPU time (in seconds) required to solve the SOS relaxation as the number $n$ of variables was increased with $d=2$. Both SOSADMM and CDCS-primal were faster than SeDuMi on these examples (note that SeDuMi's runtime reduces if a weaker termination tolerance is set, but not significantly). Also, the optimal value returned by SOSADMM was within 0.05\% of the high-accuracy value returned by SeDuMi. For all examples in Table~\ref{T:TimePOP}, the cone size $N$ is moderate (less than 300), while the number of constraints $m$ is large. SeDuMi assembles and solves an $m \times m$ linear system at each iteration, which is computationally expensive. 

\begin{table}[t]
    \centering
    \setlength{\abovecaptionskip}{0pt}
    \setlength{\belowcaptionskip}{0em}
    \renewcommand\arraystretch{1.}
    \caption{CPU time (s) to solve the SDP relaxations of \eqref{eq:ExPOP}. $N$ is the size of the PSD cone, $m$ is the number of constraints.}
    \label{T:TimePOP}
    \begin{tabular}{c c c|  c c c }
        \hline \toprule[1pt] 
        \multicolumn{3}{c|}{Dimensions} &  \multicolumn{3}{c}{CPU time (s)} \\
        \hline
         $n$   & $N$ & $m$   & SeDuMi &
         \begin{tabular}[x]{@{}c@{}}CDCS\\[-0.25em](primal)\end{tabular} &
         \begin{tabular}[x]{@{}c@{}}SOS-\\[-0.25em]ADMM\end{tabular} \\
        $2$ & 6 & 14  & 0.23  & 0.08 & 0.05  \\
        $4$ & 15 & 69  & 0.13  & 0.11 & 0.06  \\
        $6$ & 28 & 209  & 0.24  & 0.16 & 0.14  \\
        $8$ & 45 & 494  & 1.16   & 0.18 & 0.18 \\
        $10$ & 66 & 1000  & 3.17 & 0.25  & 0.39  \\
        $12$ & 91 & 1819  & 13.89  & 0.46 & 0.55  \\
        $14$ & 120 & 3059  & 54.63 & 0.79 & 0.84  \\
        $16$ & 153 & 4844 & 187.0  & 0.92 & 0.82  \\
        $18$ & 190 & 7314  & 610.2 & 2.91 & 1.92  \\
        $20$ & 231 & 10625  & 1739  & 4.93 & 2.32  \\
        \bottomrule[1pt]
        \end{tabular}
\end{table}

\subsection{Finding Lyapunov functions} \label{se:simulationLya}

\begin{table}[t]
    \centering
    \setlength{\abovecaptionskip}{0pt}
    \setlength{\belowcaptionskip}{0em}
    \renewcommand\arraystretch{1.1}
   \caption{Lyapunov functions for the system~\eqref{eq:LyapunovEx1}}
    \label{T:LyapunovEx1}
    \begin{tabular}{c | c c}
        \hline \toprule[1pt] 
        Solver &  Time (s) & Lyapunov function $V(x)$ \\
        \hline
        SeDuMi & 0.054 & $6.659x_1^2+4.628x_2^2+2.073x_3^2$   \\
        CDCS-primal & 0.21 & $7.008x_1^2+1.477x_2^2+2.172x_3^2$   \\
        SOSADMM& 0.58 & $6.699x_1^2+1.803x_2^2+2.172x_3^2$  \\
        \bottomrule[1pt]
        \end{tabular}
\end{table}

\begin{table}[t]
    \centering
    \setlength{\abovecaptionskip}{0pt}
    \setlength{\belowcaptionskip}{0em}
    \renewcommand\arraystretch{1.1}
     \caption{CPU time (s) to construct a quadratic Lyapunov function for randomly generated polynomial systems.}
    \label{T:LyapunovEx2}
    \begin{tabular}{c | c c | p{0.7cm}<{\centering} p{0.7cm}<{\centering} p{0.7cm}<{\centering}}
        \hline \toprule[1pt] 
        & \multicolumn{2}{c|}{Statistics} &  \multicolumn{3}{c}{CPU time (s)} \\
        \hline
        $n$ &
        Size of $A$&
        \begin{tabular}[x]{@{}c@{}}nonzero\\[-0.25em]density\end{tabular}
        & SeDuMi &
        \begin{tabular}[x]{@{}c@{}}CDCS\\[-0.25em](primal)\end{tabular}
        & \begin{tabular}[x]{@{}c@{}}SOS-\\[-0.25em]ADMM\end{tabular}
         \\
        10 & $1100\times 2365$ & $1.50\times 10^{-3}$ &3.3 & 7.5 & 5.3  \\
        12 & $1963\times 4407$  & $8.76\times 10^{-4}$ &11.0 & 11.8 & 7.7  \\
        14 & $3255 \times 7560$ & $5.25\times 10^{-4}$ &  49.9&21.0 & 11.2 \\
        16 & $5100 \times 12172$ & $3.13\times 10^{-4}$ &181.9 & 31.7 &16.2  \\
        18 & $7638 \times 18639$ & $2.13\times 10^{-4}$ & 574.8 & 55.0  & 24.6\\
        20 & $11025 \times 27405$ & $1.48\times 10^{-4}$ &1617.2 & 100.3 & 37.7  \\
        22 & $15433 \times 38962$ & $1.11\times 10^{-4}$ & 7442.7 & 265.9 & 65.6 \\
        25 & $24375 \times 62725$ & $6.87 \times 10^{-5}$ & * & 729.1 & 104.7 \\
        30 & $47275 \times 124620$ & $3.64 \times 10^{-5}$ & * & 3509.2 & 259.0 \\
        \bottomrule[1pt]
        \end{tabular}
        \begin{tablenotes}
      \scriptsize
      \item[\textdagger] * SeDuMi fails due to memory requirements.
    \end{tablenotes}
\end{table}

Next, we consider the problem of constructing Lyapunov functions to check local stability of polynomial/rational systems when~\eqref{eq:FindLyapunov} is replaced by SOS conditions. We used SOSTOOLS~\cite{papachristodoulou2013sostools} to generate the corresponding SDPs. 

The first system we study is
\begin{equation}\label{eq:LyapunovEx1}
    \begin{aligned}
        \dot{x}_1 &= -x^3_1 -x_1x_3^2, \\
        \dot{x}_2 &= -x_2 -x_1^2x_2, \\
        \dot{x}_3 &= -x_3 -\frac{3x_3}{x_3^2+1}+3x_1^2x_3,\\
    \end{aligned}
\end{equation}
which is demo 2 in SOSTOOLS. The system has an equilibrium at the origin, and we search for a homogeneous quadratic polynomial Lyapunov function $V(x)=ax_1^2 + bx_2^2 + cx_3^2$ to prove its global stability. The results given by SeDuMi, CDCS-primal and SOSADMM are listed in Table~\ref{T:LyapunovEx1}. For such a small system, SeDuMi was slightly faster than CDCS-primal and SOSADMM, which is expected since IPMs are well-suited for small-scale SDPs. Note that since the problem of constructing a Lyapunov functions is a feasibility problem, the solutions returned by SeDuMi, CDCS-primal and SOSADMM need not be the same (see Table~\ref{T:LyapunovEx1}).

As the last example, we consider randomly generated polynomial dynamical systems $\dot{x} = f(x)$ of degree three with a locally asymptotically stable equilibrium at the origin, and checked for local nonlinear stability in the ball $\mathcal{D} = \{x\in \mathbb{R}^n | 0.1 - \|x\|^2 \geq 0\}$ using a complete quadratic polynomial as the candidate Lyapunov function. Table~\ref{T:LyapunovEx2} summarizes the average CPU times required to search for such a Lyapunov function, when successful (note that we cannot detect infeasible problems because we only solve the primal form~\eqref{E:PrimalVectorForm}). The results clearly show that SOSADMM is faster than both SeDuMi and CDCS-primal for the largest problem instances ($n \geq 18$). Also, FOMs have much lower memory requirements, and SOSADMM can solve problems that are not accessible with IPM: SeDuMi failed due to memory issues when $n>22$. Finally, note that for the problem of finding Lyapunov functions the $m\times m$ linear system solved in SeDuMi and CDCS is not diagonal, and solving it is expensive: when $n = 30 $ it took over 150 s for CDCS just to factorize $AA^T$, which is over 50\% of the total time taken by SOSADMM to return a solution.

\section{Conclusion} \label{se:conclusion}

In this paper, we proposed an efficient ADMM algorithm to exploit the row-sparsity of SDPs that arise from SOS programming, which are implemented in SOSADMM. The subproblems of our algorithm consist of one conic projection and multiple quadratic programs with closed-form solutions, which can be computed efficiently and---most importantly--- do not require any matrix inversion.

Our numerical experiments on random unconstrained polynomial optimization and on Lyapunov stability analysis of polynomial/rational systems demonstrate that our method can provide speed-ups compared to the interior-point solver SeDuMi and the first-order solver CDCS. One major drawback of our method is the inability to detect infeasibility; future work will try to exploit the sparsity of SDPs from SOS relaxations in a homogeneous self-dual embedding formulation similar to that of~\cite{ODonoghue2016,zheng2016fast}.

\bibliographystyle{IEEEtran}
\bibliography{Reference}

\end{document}